\documentclass{article}
\usepackage[a4paper, margin=1in]{geometry}
\usepackage{amsmath,amsfonts,amssymb,amsthm}
\usepackage{braket}
\usepackage{graphicx}
\usepackage[dvipsnames]{xcolor}
\usepackage{subcaption}
\usepackage{float}
\usepackage{tabularx}
\usepackage{array}
\usepackage{booktabs}
\usepackage{multirow}
\usepackage{tikz}
\usepackage{siunitx}
\usepackage[table]{xcolor}
\usepackage{makecell}
\usepackage[linesnumbered,ruled,vlined]{algorithm2e}
\SetKwInput{Input}{Input}
\SetKwInput{Output}{Output}
\usepackage[normalem]{ulem}
\usepackage[font=small]{caption}
\usepackage{url}

\newcommand{\sA}{\textcolor{RedViolet}{S_\text{A}}}
\newcommand{\sP}{\textcolor{MidnightBlue}{S_\text{P}}}
\newcommand{\scat}{\textcolor{PineGreen}{S_\text{cat}}}
\newcommand{\scatA}{\textcolor{Sepia}{S_\text{catA}}}

\newcommand{\xA}{\textcolor{RedViolet}{x_\text{A}}}
\newcommand{\xP}{\textcolor{MidnightBlue}{x_\text{P}}}
\newcommand{\xcat}{\textcolor{PineGreen}{x_\text{cat}}}
\newcommand{\xcatA}{\textcolor{Sepia}{x_\text{catA}}}

\newcommand{\dxA}{\textcolor{RedViolet}{\dot{x_\text{A}}}}
\newcommand{\dxP}{\textcolor{MidnightBlue}{\dot{x_\text{P}}}}
\newcommand{\dxcat}{\textcolor{PineGreen}{\dot{x_\text{cat}}}}
\newcommand{\dxcatA}{\textcolor{Sepia}{\dot{x_\text{catA}}}}

\title{Integral chemical reaction neural networks}
\author{Abraham Reyes-Vel\'{a}zquez\thanks{Department of Mathematics, The University of Manchester, Oxford Road, Manchester, M13\,9PL, United Kingdom, \texttt{abrahamrafael.reyesvelazquez@manchester.ac.uk}} \and Stefan G\"{u}ttel\thanks{Department of Mathematics, The University of Manchester, Oxford Road, Manchester, M13\,9PL, United Kingdom, \texttt{stefan.guettel@manchester.ac.uk}}}
\date{\today}

\begin{document}

\maketitle

\begin{abstract}
    Discovering the structure and kinetics of chemical reaction networks (CRNs) from time-series concentration data is a fundamental challenge in chemical kinetics, with existing approaches relying on prior mechanistic assumptions or suffering from high computational cost and noise sensitivity. In this work, we present integral chemical reaction neural networks (iCRNNs), a framework that combines the interpretable, physics-constrained, architecture of chemical reaction neural networks (CRNNs) with an integral collocation formulation of the governing dynamics. Rather than solving ODEs forward in time at each training step, we approximate the integral of the learned rate functions directly using numerical quadrature. This yields an entirely algebraic forward pass consisting only of matrix operations, eliminating the repeated adaptive ODE solves of standard CRNN and producing smoother, more predictable training. We provide a recovery error analysis characterising the structural sources of ill-conditioning (conservation laws, reaction reversibility, and shared reactant sets) that limit network identifiability. On two benchmark CRNs, iCRNN trains in roughly half the wall-clock time of the baseline CRNN method on a four-species mechanism and between three and four times faster on a five-species mechanism, while completing every training run reliably and attaining comparable or lower loss.
\end{abstract}

\section{Introduction} \label{sec:intro}
\noindent Identifying the chemical mechanism occurring within a reservoir, given only time-series data of its chemical species concentrations, remains a central challenge in chemical reaction network theory. Most chemical system identification methods focus on automating the extraction of the rate laws, the governing ordinary differential equations (ODEs) that dictate how species concentrations change over time. This includes classical kinetic analysis tools such as the integral method~\cite{Espenson1995}, the initial rates method~\cite{Espenson1995}, the Delplot technique~\cite{bhore1990delplot}, reaction progress kinetic analysis (RPKA)~\cite{blackmond2005reaction}, and variable time normalisation analysis (VTNA)~\cite{bures2016variable}, as well as recent data-driven frameworks such as two-step collocation methods, including sparse identification of nonlinear dynamics (SINDy)~\cite{brunton2016discovering}, and graph reconstruction via additive differential equations (GRADE)~\cite{chen2017network}. However, mapping these abstract empirical rate laws back to a physically meaningful stoichiometric graph is still not automated, and is thus still performed heuristically through the intuition of an experienced chemist.

Among the frameworks that explicitly attempt to recover valid stoichiometric mechanisms, optimisation methods and machine learning architectures have shown promise. For example, neural ordinary differential equations (NODEs) provide a continuous-time framework for learning complex dynamical systems directly from their trajectories without requiring a fully specified a priori physical model~\cite{chen2018neural}. Within CRN theory, this approach took form in the chemical reaction neural network (CRNN) architecture proposed by Ji \& Deng~\cite{ji2021autonomous}. This architecture hardcodes the mass-action law and can thus produce models with a direct stoichiometric interpretation. However, there are instrinsic limitations: several CRN structures can obscure the model discovery process, and since CRNN training constitutes a continuous relaxation of the discrete CRN discovery problem, we may obtain only a good trajectory predictor. Moreover, the performance of this approach depends strongly on two key components: the numerical method used to solve an initial value problem (IVP) during each forward pass and the sensitivity analysis method used during the backward pass. Different combinations of these methods are better suited to different problems~\cite{ma2021comparison, lakrisenko2024benchmarking}. In summary, CRN discovery via CRNN training faces three challenges, two numerical and one structural:

\begin{enumerate}
    \item if the CRN exhibits a large disparity in reaction rates, with both \emph{slow} and \emph{fast} reactions, the resulting trajectories can be stiff, requiring computationally expensive implicit solvers at each forward pass;

    \item if the CRN contains conservation laws between species concentrations, its Jacobian can become ill-conditioned~\cite{oellerich2021biological}, requiring special care in sensitivity analysis~\cite{lakrisenko2024benchmarking};
    
    \item if the CRN contains certain reaction structures, such as reversible reactions or multiple reactions sharing the same reactant, CRNN training may fail to recover the correct stoichiometric structure.
\end{enumerate}

To mitigate the numerical difficulties associated with the first two challenges, we propose iCRNN, a collocation method~\cite{varah1982spline} in which the forward-pass integration is realised algebraically through a matrix product involving an integration matrix $\mathbf{J}$. This matrix is constructed by numerical quadrature of cubic spline interpolants to the training data. Because the forward pass is an explicit algebraic discretisation, its derivatives with respect to the model parameters can be evaluated directly by automatic differentiation, without differentiating through an adaptive ODE solver. Furthermore, it has been shown~\cite{RGLL26} that this integral approximation is more robust to measurement noise than methods based on differentiation. The third challenge is structural rather than numerical; we characterise its sources and consequences, showing in particular that shared reactants can render the ground-truth CRN unidentifiable while the model still fits the data. Specifically, this work makes the following contributions:
\begin{enumerate}
    \item We introduce iCRNN, a reformulation of CRNN in algebraic form. Rather than repeatedly solving the learned ODE during training, we propose a purely algebraic forward pass consisting entirely of matrix operations. iCRNN evaluates the learned reaction rates at the observed concentration time points and applies a fixed spline-based collocation operator to approximate their time integrals. This eliminates the need for repeated adaptive ODE solves during training and yields a smooth, predictable loss history.

    \item We provide a recovery error analysis that decomposes the parameter recovery error into a conditioning factor determined by the network structure and a data-error factor determined by the noise level and the accuracy of the collocation approximation. We identify three structural sources of ill-conditioning (conservation laws, reaction reversibility, and shared reactants) and show that shared reactants are qualitatively distinct, since they can render the true network topology unidentifiable while the learned model still fits the data.

    \item We validate the approach on two benchmark CRN systems, demonstrating a runtime speedup, completion of training across all runs, and loss performance comparable to or better than baseline CRNN.
\end{enumerate}

This paper is organised as follows. Section~\ref{sec:crn} provides an overview of chemical reaction network theory and introduces the notation used throughout. Section~\ref{sec:crnn} describes the CRNN architecture developed by Ji \& Deng~\cite{ji2021autonomous}. Section~\ref{sec:icrnn} presents our iCRNN framework and details the training algorithm. Section~\ref{sec:error} provides a recovery error analysis characterising structural sources of ill-conditioning. Section~\ref{sec:experiments} presents numerical experiments validating the iCRNN framework. Finally, Section~\ref{sec:conclusions} concludes with a discussion of limitations and future work.

\section{Chemical reaction networks}\label{sec:crn}
Chemical reaction networks (CRNs) are many-body dynamical systems that describe the interactions and transformations of discrete populations of particles. These systems consist of a collection of particles contained in a reservoir. These particles can be of different types, each moving randomly within the reservoir. Transformations of particles of one species into particles of another species are modelled through collisions: when a set of reacting particles collides in the reservoir, a reaction is triggered, thus producing a set of product particles. Because of the generality of this framework, CRNs are widely used in chemistry, biology and epidemiology to represent population dynamics driven by discrete events. Formally, a CRN is defined by a triple of sets \((\mathcal{S},\mathcal{Q},\mathcal{R})\) such that:
\begin{itemize}
    \item \(\mathcal{S}=\{S_{\alpha}\}_{\alpha=1}^{M}\) is the set of species, with an element for each different kind of particle we can have in our system;
    \item \(\mathcal{Q} = \{q_{i}\}_{i=1}^{N}\) is the set of complexes, whose elements are the different ways in which the particles can interact to perform a chemical reaction. Mathematically, the set \(\mathcal{Q}\) contains \emph{all possible} linear combinations \(q_{i} = \sum_{\alpha=1}^{M}Q_{\alpha,i} S_{\alpha}\) with integer weights \(Q_{\alpha,i}\geq 0\) such that \(\sum_{\alpha=1}^{M} Q_{\alpha,i} \leq p\) for some total degree \(p\) and with at least one $Q_{\alpha,i}$ nonzero for each $i$. The degree $p$ corresponds to the highest reaction order in the network, i.e., the size of the largest reactant complex (in number of particles), and we have the relation
    \[
        N = \binom{M+p}{p}-1;
    \]
    \item \(\mathcal{R}\) is the set of directed reactions \((q_i,q_j)\) where each ordered pair represents the transformation of complex~\(q_{i}\) into complex~\(q_{j}\).
\end{itemize}

CRNs are inherently stochastic systems driven by discrete events, and the mathematical object of primary interest is the time-evolution of the species population vector \(\mathbf{n}(t) = (n_{1}(t),\ldots,n_M(t)) \in \mathbb{Z}_{\geq 0}^{M}\). Under appropriate conditions, namely, that the particles are homogeneously mixed in the reservoir and the total population is sufficiently large, the stochastic description can be approximated by a continuous, deterministic one. At this thermodynamic limit, the dynamics can be described by the evolution of the species concentration vector $ \mathbf{x}(t) \in \mathbb{R}^{M}$. Under this deterministic regime, CRNs satisfy mass-action kinetics where the rate of reaction $(q_i, q_j) \in \mathcal{R}$ is
\[
    r_{(i,j)}(\mathbf{x}(t)) = k_{i,j} \prod_{\alpha=1}^{M} [x_\alpha(t)]^{Q_{\alpha,i}},
\]
where $k_{i,j} \geq 0$ is the rate constant of reaction $(q_i,q_j)$, and $Q_{\alpha,i}$ is the stoichiometric coefficient of species $S_\alpha$ in the reactant complex $q_i$. The contribution of this reaction to the time-evolution of species $\alpha$'s concentration is the kinetic term
\[
    \kappa_{\alpha,(i,j)}(\mathbf{x}(t)) 
    = \bigl[Q_{\alpha,j} - Q_{\alpha,i}\bigr]\, r_{(i,j)}(\mathbf{x}(t)).
\]
The dynamics of the CRN are then governed by the (generally nonlinear) ODE system
\[
    \frac{d}{dt} x_\alpha(t) = \sum_{(i,j) \in \mathcal{R}} 
    \kappa_{\alpha,(i,j)}(\mathbf{x}(t)),
    \quad \alpha \in \{1,\ldots,M\}.
\]
To write the full system compactly, let $\mathbf{Q} = [Q_{\alpha,i}] \in \mathbb{Z}_{\geq 0}^{M \times N}$ be the complex stoichiometry matrix, and let $\mathbf{d}(\mathbf{x}(t)) \in \mathbb{R}^N$ be the mass-action monomial dictionary with entries
\[
    d_i = \prod_{\alpha=1}^{M} [x_\alpha(t)]^{Q_{\alpha,i}}, 
    \quad i \in \{1,2,\ldots,N\}.
\]
We remark that $\mathbf{d}(\mathbf{x}(t))$ has an entry per each complex in $\mathcal{Q}$. Now construct the Kirchhoff matrix $\mathbf{K} \in \mathbb{R}^{N \times N}$ with off-diagonal entries $k_{i,j} \geq 0$ and diagonal entries $-\sum_{j} k_{i,j}$, so that every column sums to zero. Then the ODE system of the network, in graph-explicit form, reads
\begin{equation}\label{eq:crn-ode}
    \dot{\mathbf{x}}(t) = \mathbf{Q}\mathbf{K}\,\mathbf{d}(\mathbf{x}(t)), 
    \qquad \mathbf{x}(0) \text{ given.}
\end{equation}
The coefficient matrix $\mathbf{C} = \mathbf{Q}\mathbf{K} \in \mathbb{R}^{M \times N}$ is the target of dictionary-based recovery methods such as SINDy~\cite{brunton2016discovering}. This parametrisation, however, conflates the network graph structure with the reaction kinetics, obscuring chemical interpretability.

To expose the reaction structure more explicitly, we decompose the Kirchhoff matrix as $\mathbf{K} = \boldsymbol{\Lambda}^{(1)} \boldsymbol{\Lambda}^{(2)}$, where $\boldsymbol{\Lambda}^{(1)} \in \mathbb{Z}^{N \times R}$ encodes network connectivity and $\boldsymbol{\Lambda}^{(2)} \in \mathbb{R}^{R \times N}$ encodes rate constants, with $R = |\mathcal{R}|$ the number of reactions. Concretely,
\[
    \Lambda^{(1)}_{\ell,(i,j)} = \delta_{\ell,j}-\delta_{\ell,i}, 
    \quad \ell \in \{1,2,\ldots,N \},
    \quad (i,j) \in \{1,2,\ldots,R \}
\]
(i.e., $-1$ if $q_l$ is the reactant complex of the reaction, $+1$ if it is the product complex, and $0$ otherwise), and
\[
    \Lambda^{(2)}_{(i,j),\ell} = \delta_{\ell,i}\, k_{i,j}.
\]
The product $\mathbf{Q}^{(\Delta)} = \mathbf{Q}\boldsymbol{\Lambda}^{(1)} \in \mathbb{Z}_{}^{M \times R}$ gives the reaction stoichiometry matrix, encoding the net change in each species per reaction, with entries
\[
    Q^{(\Delta)}_{\alpha,(i,j)} = Q_{\alpha,j} - Q_{\alpha,i}.
\]
Consider the product stoichiometry matrix $\mathbf{Q}^{(\mathrm{Pr})} = [Q_{\alpha,j}] \in \mathbb{Z}_{\geq 0}^{M \times R}$ and the reactant stoichiometry matrix $\mathbf{Q}^{(\mathrm{Re})} = [Q_{\alpha,i}] \in \mathbb{Z}_{\geq 0}^{M \times R}$. We these we can also write the reaction stoichiometry matrix as $\mathbf{Q}^{(\Delta)} = \mathbf{Q}^{(\mathrm{Pr})} - \mathbf{Q}^{(\mathrm{Re})}$. 

The product $\mathbf{r}(\mathbf{x}(t)) = \boldsymbol{\Lambda}^{(2)} \mathbf{d}(\mathbf{x}(t)) \in \mathbb{R}^R$ yields the rates vector with entries $r_{(i,j)}(\mathbf{x})$, so that the dynamics~\eqref{eq:crn-ode} take the reaction-explicit form
\begin{equation}\label{eq:crn-rates}
    \dot{\mathbf{x}}(t) = \mathbf{Q}^{(\Delta)}\,\mathbf{r}(\mathbf{x}(t)), 
    \qquad \mathbf{x}(0) \text{ given.}
\end{equation}
The reaction-explicit representation~\eqref{eq:crn-rates} separates network structure (encoded in $\mathbf{Q}^{(\Delta)}$) from kinetics (encoded in $\mathbf{r}$). 

In Table~\ref{tab:matrix_interpretation} we summarise the chemical interpretations of the matrices and vectors introduced here.

\begin{table}[H]
    \centering
    \caption{Key matrices and vectors in chemical reaction network theory and their physical interpretation}
    \label{tab:matrix_interpretation}
    
    \arrayrulecolor{gray!60}
    
    \begin{tabular}{|c|p{10cm}|}
        \hline
        \rowcolor{gray!15}
        \textbf{Matrix or vector} &
        \textbf{Interpretation} \\
        \hline
        
        $\mathbf{Q} \in \mathbb{Z}_{\geq 0}^{M \times N}$&
        Complex stoichiometry: $Q_{\alpha,i}$ is the number of units of species $\alpha$ in complex $q_i$. \\
        \hline
        
        $\mathbf{Q}^{(\mathrm{Re})} \in \mathbb{Z}_{\geq 0}^{M \times R}$&
        Reactant stoichiometry: $Q^{(\mathrm{Re})}_{\alpha,(i,j)}$ is the number of units of species $\alpha$ in reactant complex $q_i$. \\
        \hline
        
        $\mathbf{Q}^{(\mathrm{Pr})} \in \mathbb{Z}_{\geq 0}^{M \times R}$&
        Product stoichiometry: $Q^{(\mathrm{Pr})}_{\alpha,(i,j)}$ is the number of units of species $\alpha$ in product complex $q_j$. \\
        \hline
        
        $\mathbf{Q}^{(\Delta)} \in \mathbb{Z}_{}^{M \times R}$&
        Reaction stoichiometry: $Q^{(\Delta)}_{\alpha,(i,j)}$ is the net change in species $\alpha$ due to reaction $(q_i,q_j)$. \\
        \hline
        
        $\mathbf{K} \in \mathbb{R}^{N \times N}$
        &
        Kirchhoff matrix: encodes reaction rates between complexes. For $i \neq j$, the off-diagonal element $K_{i,j}=k_{i,j}$ is the kinetic constant of reaction $(q_i,q_j)$, while the diagonal element $K_{i,i}=-\sum_{j\neq i}k_{i,j}$ is the negative sum of the kinetic constants of all reactions originating from complex $q_i$. \\
        \hline
        
        $\boldsymbol{\Lambda}^{(1)} \in \mathbb{Z}_{}^{N \times R}$&
        Network connectivity matrix: $\Lambda^{(1)}_{\ell,(i,j)}$ indicates whether complex $\ell$ is the reactant ($-1$), product ($+1$), or neither ($0$) in reaction $(q_i,q_j)$. \\
        \hline
        
        $\boldsymbol{\Lambda}^{(2)} \in \mathbb{R}^{R \times N}$
        &
        Rate constant matrix: $(\Lambda^{(2)})_{(i,j),l}$ encodes the kinetic constants, placing $k_{i,j}$ at the row corresponding to reaction $(q_i,q_j)$ and the column corresponding to reactant complex $q_i$. \\
        \hline
        
        $\mathbf{x} \in \mathbb{R}^M$
        &
        Concentration vector. \\
        \hline
        
        $\mathbf{d}(\mathbf{x}) \in \mathbb{R}^N$
        &
        Mass-action monomial dictionary vector: $$d_i(\mathbf{x})=\prod_{\alpha=1}^M x_\alpha^{Q_{\alpha,i}}.$$\\
        \hline
        
        $\mathbf{r}(\mathbf{x}) \in \mathbb{R}^R$
        &
        Mass-action reaction-rate vector: $$r_{(i,j)}(\mathbf{x})=k_{i,j}\prod_{\alpha=1}^M x_\alpha^{Q_{\alpha,i}}.$$\\
        \hline
        
    \end{tabular}
\end{table}

\section{Chemical reaction neural networks}\label{sec:crnn}
The chemical reaction neural network (CRNN)~\cite{ji2021autonomous} discovers the reaction network directly from time-series concentration data, without prior knowledge of which reactions are present or what their rate constants are.

\subsection{Model parametrisation}
The CRNN assumes that the system takes the structural form~\eqref{eq:crn-rates} and parametrises each reaction rate via the log-exp transform~\cite{ji2021autonomous}:
\begin{align*}
    r_{(i,j)}(\mathbf{x}(t)) 
    &= \exp \Bigg\{ \log \Bigg[k_{i,j} \prod_{\alpha=1}^{M} \Big(x_{\alpha}(t)\Big)^{Q_{\alpha,i}} \Bigg] \Bigg\} \\
    &= \exp \Bigg( \sum_{\alpha=1}^{M} Q_{\alpha,i}\, \log(x_{\alpha}(t)) 
    + \log(k_{i,j}) \Bigg) \\[7pt]
    &= \exp \Bigg( \Big(\mathbf{Q}^{(\mathrm{Re})}_{:,(i,j)} \Big)^{\mathrm{T}} 
    \log(\mathbf{x}(t)) + b_{(i,j)}\Bigg),
\end{align*}
where $b_{(i,j)} = \log( k_{i,j})$. Stacking all reaction rates, the full reparametrised rates vector reads
\[
    \mathbf{r}(\mathbf{x}(t)) 
    = \exp\Bigg(\Big(\mathbf{Q}^{(\mathrm{Re})} \Big)^{\mathrm{T}} 
    \log(\mathbf{x}(t)) + \mathbf{b}\Bigg),
\]
where $\mathbf{b} \in \mathbb{R}^R$ collects the log-rate constants. With this reparametrisation, the Picard integral form of the reaction-explicit ODE system~\eqref{eq:crn-rates} becomes the neural ODE~\cite{chen2018neural}
\begin{equation} \label{eq:NODE}
    \mathbf{x}(t) = \mathbf{x}(0) + \int_{0}^{t} 
    f(\mathbf{x}(\tau), \tau,\boldsymbol{\theta})\, d\tau, 
    \quad \boldsymbol{\theta}=\left[\mathrm{vec}\left(\mathbf{Q}^{(\Delta)}\right), \mathrm{vec}\left(\mathbf{Q}^{(\mathrm{Re})}\right),\, \mathbf{b}\right]^{\mathrm{T}} 
    \in \mathbb{R}^{R(2M+1)},
\end{equation}
where the notation $\mathrm{vec}(\cdot)$ denotes a flattened vector constructed from the given matrix, and
\[
    f(\mathbf{x}(\tau), \tau,\boldsymbol{\theta}) = \mathbf{Q}^{(\Delta)}
    \exp\Bigg(\Big(\mathbf{Q}^{(\mathrm{Re})} \Big)^{\mathrm{T}} 
    \log(\mathbf{x}(\tau)) + \mathbf{b}\Bigg).
\]
The size of the parameter vector $\boldsymbol{\theta}$ depends on the number of species $M$ and the number of reactions $R$ in the CRN. While $M$ is determined by the number of features in the training data, the choice of $R$ is left as a hyperparameter that might require tuning if we do not have any a priori knowledge about the CRN.

In the CRN formalism, the entries of $\mathbf{Q}^{(\mathrm{Re})}$ and $\mathbf{Q}^{(\Delta)}$ are integer-valued and satisfy structural constraints inherited from the complex graph $\mathbf{K}$. The CRNN replaces this discrete parameter space by a continuous relaxation in which these matrices entries are optimised over $\mathbb{R}$. Consequently, a trained CRNN might not correspond exactly to a valid chemical reaction network, even when it accurately reproduces the observed dynamics.

\subsection{Training}
Suppose that, for a given initial condition $\mathbf{x}^{(w)}(0)$, we observe trajectories of all $M$ species of a CRN at time points $0 = t_0 < t_1 < \cdots < t_n$. The noisy measurements are given by
\[
    \bar{x}_{\alpha,\ell}^{(w)} = x_\alpha^{(w)}(t_\ell) + \xi_{\alpha,\ell}^{(w)},
\]
where the measurement noise $\xi_{\alpha,\ell}^{(w)}$ is assumed to be i.i.d. and globally bounded, i.e., there exists $\varepsilon>0$ such that
\[
    |\xi_{\alpha,\ell}^{(w)}|\leq \varepsilon
\]
for all species $\alpha$, time points $t_\ell$, and experiments $w$.

The measurements from experiment $w$ are collected into the data matrix $\mathbf{X}^{(w)} = [\bar{x}_{\alpha,\ell}^{(w)}] \in \mathbb{R}^{M \times (n+1)}$. We consider $W$ experiments in total and divide them into $W_{\mathrm{train}}$ training experiments and $W_{\mathrm{val}}$ validation experiments, with
\[
    W=W_{\mathrm{train}}+W_{\mathrm{val}}.
\]

Training proceeds by sequentially processing the training experiments:

Given the current parameter vector $\boldsymbol{\theta}_\ell$ and initial conditions $\mathbf{x}^{(w)}(0)$ of the $w$-th training experiment, we generate a prediction $\widehat{\mathbf{X}}_\ell^{(w)}$. We compare this prediction with the observed data $\mathbf{X}^{(w)}$ through a loss function 
\[
    \mathcal{L}(\boldsymbol{\theta}_\ell, \ \widehat{\mathbf{X}}^{(w)}_{\ell}, \ \mathbf{X}^{(w)}) = \mathcal{L}^{(w)}(\boldsymbol{\theta}_{\ell}).
\]

Two loss functions are used in this work. For the error analysis of Section~\ref{sec:error}, we require a loss that is twice differentiable and admits a closed-form Hessian. We therefore adopt the mean squared error (MSE)
\begin{equation} \label{eq:mse-loss}
    \mathcal{L}^{(w)}_{\mathrm{MSE}}(\boldsymbol{\theta}_{\ell}) = \frac{1}{M(n+1)} \left\| \mathbf{X}^{(w)} - \widehat{\mathbf{X}}^{(w)}_{\ell} \right\|_F^2.
\end{equation}
For training in Section~\ref{sec:experiments}, we follow \cite{ji2021autonomous} and use the mean absolute error (MAE)
\begin{equation} \label{eq:training-loss}
    \mathcal{L}^{(w)}_{\mathrm{MAE}}(\boldsymbol{\theta}_{\ell}) = \frac{1}{M(n+1)} \left\| \mathbf{X}^{(w)} - \widehat{\mathbf{X}}^{(w)}_{\ell} \right\|_1,
\end{equation}
where the entry-wise $\ell^1$ norm is defined by
\[
    \|\mathbf{Z}\|_1 = \sum_{i,j}|Z_{i,j}|,
\]
as well as the MSE~\eqref{eq:mse-loss}. The gradient
\[
    \nabla_{\boldsymbol{\theta}_\ell}\mathcal{L}^{(w)}(\boldsymbol{\theta}_\ell)
\]
is computed by automatic differentiation through the forward pass, and the parameters are subsequently updated using the AdamW optimiser~\cite{loshchilov2018decoupled}. Each training experiment corresponds to one parameter update. Processing all $W_{\mathrm{train}}$ training experiments completes an epoch. At the end of each epoch, the updated parameters are used to generate predictions for all training and validation experiments, and the mean training and validation losses are evaluated over their respective experiments. The validation experiments are not used to update the model parameters.

Training the neural ODE in~\eqref{eq:NODE} by updating the parameters $\boldsymbol{\theta}$ yields a model with an explicitly parametrised reaction structure, allowing us to infer:
\begin{itemize}
    \item the candidate reaction structure (through the sparsity pattern of $\mathbf{Q}^{(\Delta)}$);
    \item the stoichiometry of each candidate reactant (through the values of $\mathbf{Q}^{(\mathrm{Re})}$);
    \item the rate constants of each candidate reactant (through the values of $\mathbf{k} = \exp(\mathbf{b})$).
\end{itemize}

In standard CRNN, the integral in~\eqref{eq:NODE} is evaluated numerically via an ODE solver (e.g., Runge--Kutta methods such as RK4 or  Tsit5~\cite{tsitouras2011runge}), to produce the predicted trajectories $\widehat{\mathbf{X}}_\ell$ at the required time points. Gradients of the training loss $\mathcal{L}^{(w)}(\boldsymbol{\theta}_\ell)$ with respect to $\boldsymbol{\theta}_\ell$ are then computed either by back propagating through the solver operations or via the adjoint sensitivity method. This approach carries several practical drawbacks. Each forward pass requires solving a full ODE, which is expensive when the integration must be repeated over many training iterations and experiments. Stiff reaction networks may demand very small step sizes, and adaptive solvers introduce further dependence on tolerance hyperparameters and can produce memory-intensive gradient computations via the adjoint method.

\section{iCRNN: integral formulation of CRNNs}\label{sec:icrnn}
We propose an alternative, algebraic, forward pass that avoids ODE solvers entirely. Given current parameters $\boldsymbol{\theta}$, the rates matrix is
\[
    \mathbf{R}(\mathbf{X}) 
    = \exp \bigg(\Big(\mathbf{Q}^{(\mathrm{Re})}\Big)^{\mathrm{T}}\log(\mathbf{X}) 
    + \mathbf{b}\,\mathbf{1}_{n+1}^\mathrm{T}\bigg)
    \in \mathbb{R}^{R \times (n+1)},
\]
where $\log(\mathbf{X})$ is applied entry-wise, and the data-level dynamical system then reads
\begin{equation*}\label{eq:crn-data}
    \dot{\mathbf{X}} = \mathbf{Q}^{(\Delta)}\,\mathbf{R}(\mathbf{X}),
\end{equation*}
where $\dot{\mathbf{X}} \in \mathbb{R}^{M \times (n+1)}$ is the matrix of derivatives evaluated at the observed concentration states.

Let $\mathbf{J} \in \mathbb{R}^{(n+1)\times(n+1)}$ be a cubic spline integration matrix \cite{de1978practical} with entries defined as
\[
    J_{j,k} = \int_{t_0}^{t_j} s_k(\tau)\,d\tau,
\]
where $s_k(\tau)$ is the cardinal cubic spline basis function satisfying $s_k(t_i) = \delta_{ik}$. This matrix is constructed once from the observed time points and reused throughout training. The predicted trajectory matrix is then
\begin{equation}\label{eq:icrnn-forward}
    \widehat{\mathbf{X}} = \mathbf{X}_0 + \mathbf{Q}^{(\Delta)}\,\mathbf{R}(\mathbf{X})\,\mathbf{J},
\end{equation}
where $\mathbf{X}_0 = \mathbf{x}(t_0)\mathbf{1}_{n+1}^\mathrm{T} \in \mathbb{R}^{M \times (n+1)}$ broadcasts the initial condition across all time steps. The entire forward pass~\eqref{eq:icrnn-forward} consists of matrix multiplications and element-wise operations, with no ODE solve required. This offers several advantages over the CRNN baseline:
\begin{itemize}
    \item \textbf{Computational efficiency.} The rates matrix $\mathbf{R}(\mathbf{X})$ is evaluated once per forward pass, directly from the discrete observations, and integrated via a single matrix multiplication. In contrast, the neural ODE approach requires evaluating $\mathbf{r}(\mathbf{x}(\tau))$ at many intermediate time points per step.

    \item \textbf{Stable gradients.} The forward pass is a fixed algebraic expression, making the loss analytically differentiable via standard automatic differentiation, with no dependence on solver tolerance or stiffness.

    \item \textbf{Noise robustness.} The integral formulation naturally averages over the observed trajectory, rather than propagating errors forward from initial conditions.
\end{itemize}

Applying~\eqref{eq:icrnn-forward} to each training experiment yields the predicted trajectories $\widehat{\mathbf{X}}^{(w)}$, which are compared to the observed trajectories $\mathbf{X}^{(w)}$ through the training loss $\mathcal{L}^{(w)}(\boldsymbol{\theta}_\ell)$. During each epoch, the training experiments are processed sequentially, with the model parameters updated after each experiment. Thus, each epoch consists of $W_{\mathrm{train}}$ parameter updates. The updated parameters are subsequently used to evaluate the training and validation losses. The complete iCRNN training procedure is summarised in Algorithm~\ref{alg:icrnn-training}.

\begin{algorithm}
    \caption{iCRNN training} \label{alg:icrnn-training}
    \Input{Multi-experiment concentration time-series matrices $\{\mathbf{X}^{(w)}\}_{w=1}^{W}$, $\mathbf{X}^{(w)} \in \mathbb{R}^{M \times (n+1)}$ divided into training experiments $\{\mathbf{X}^{(w)}\}_{w=1}^{W_{\mathrm{train}}}$, and validation experiments $\{\mathbf{X}^{(w)}\}_{w=W_{\mathrm{train}}+1}^{W}$; integration matrix $\mathbf{J}$; number of epochs $N_{\mathrm{epoch}}$; loss  $\mathcal{L}^{(w)}(\boldsymbol{\theta}_\ell)$.}
    \vspace{2pt}
    \Output{Learned model $\boldsymbol{\theta}=\left[\mathrm{vec}\left(\mathbf{Q}^{(\Delta)}\right), \mathrm{vec}\left(\mathbf{Q}^{(\mathrm{Re})}\right),\, \mathbf{b}\right]^{\mathrm{T}} \in\mathbb{R}^{R(2M+1)}$.}
    \vspace{2pt}
    Initialise $\boldsymbol{\theta}$ with random entries.
    
    \vspace{2pt}
    \For{$e=1,\ldots,N_{\mathrm{epoch}}$}{
        \vspace{2pt}
        \For{$w=1,\ldots,W_{\mathrm{train}}$}{
            \vspace{2pt}
            Compute the rates matrix
            \[
                \mathbf{R}^{(w)} = \exp\left(\left(\mathbf{Q}^{(\mathrm{Re})}\right)^{\mathrm{T}} \log\left(\mathbf{X}^{(w)}\right) + \mathbf{b}\mathbf{1}_{n+1}^{\mathrm{T}} \right).
            \]
    
            Compute the predicted trajectory
            \[
                \widehat{\mathbf{X}}^{(w)} = \mathbf{X}^{(w)}_0 + \mathbf{Q}^{(\Delta)} \mathbf{R}^{(w)} \mathbf{J}.
            \]
    
            Evaluate the experiment loss $\mathcal{L}^{(w)}(\boldsymbol{\theta}_\ell)$.
    
            Compute $\nabla_{\boldsymbol{\theta}_\ell}\mathcal{L}^{(w)}(\boldsymbol{\theta}_\ell)$using automatic differentiation.
    
            Update $\boldsymbol{\theta}$.
        }
        
        Compute the loss for each training and validation experiment using the updated parameters.
        
        Compute and store the mean training and validation losses.
    }
    
    \Return{$\boldsymbol{\theta}$}

\end{algorithm}

\section{Error analysis} \label{sec:error}

The recovery of a chemical model from concentration measurements is affected by two distinct kinds of error: numerical and structural. First, the measured trajectories in our dataset may be noisy, and the integral operator approximation might fail to approach the exact time integral. Second, even in the absence of these numerical sources of error, the mapping from the iCRNN parameters to a valid chemical reaction network graph topology might be obscured due to linear dependencies between the effects of different reactions. In this section, we characterise these sources of error. We provide a bound on the iCRNN recovery error, and we identify and analyse the structural conditions under which the recovery problem becomes ill-conditioned.

\subsection{Parameter recovery error}
Let \( \boldsymbol{\theta}^{*} \in \mathbb{R}^{d}, \ d=R(2M+1), \) be a ground-truth parameter vector we want to recover. Let $\hat{\boldsymbol{\theta}} \in \mathbb{R}^{d}$ be our iCRNN-learned estimation. We want to quantify the parameter recovery error
\[
    \| \boldsymbol{\theta}^{*} - \hat{\boldsymbol{\theta}} \|_{2}.
\]

For a given experiment, we define the exact forward integral map
\[
    \mathcal{I}(\boldsymbol{\theta}) =
    \mathbf{X}_0 +  \int_0^{t} f(\mathbf{x}(\tau),\tau,\boldsymbol{\theta}) \,d\tau,
\]
and the iCRNN collocation forward integral map
\[
    \widehat{\mathcal{I}}_{\mathbf J}(\boldsymbol{\theta}) = \mathbf{X}_0 + \mathbf{Q}^{(\Delta)} \mathbf{R}(\mathbf X) \mathbf J.
\]
Their difference gives the collocation error at any parameter $\boldsymbol{\theta}$
\[
    \mathbf{E}_{\text{coll}}(\boldsymbol{\theta}) = \widehat{\mathcal{I}}_{\mathbf{J}}(\boldsymbol{\theta}) - \mathcal{I}(\boldsymbol{\theta}).
\]

From the integral approximation result established in \cite[Corr.~1]{RGLL26}, this error can be bounded in terms of the spline approximation error and the temporal discretisation,
\begin{equation} \label{eq:coll-bound}
    \|\mathbf{E}_{\text{coll}}(\boldsymbol{\theta})\|_F \le \sqrt{M(n+1)} \left( \frac{\kappa_{\mathrm{int}}}{n^4} + \varepsilon n^{1/2} \right),
\end{equation}
where \(\kappa_{\mathrm{int}}\) depends on the fourth derivatives of the dictionary functions encoded on $\mathbf{R}(\mathbf{X})$, and \(\varepsilon\) is a global bound on the noise level.

Let $\mathbf{X} = \mathcal{I}(\boldsymbol{\theta}^{*}) + \mathbf{E}_{\text{noise}}$ be the noisy measurements with $\mathcal{I}(\boldsymbol{\theta}^{*})$ the ground-truth trajectory and $\mathbf{E}_{\text{noise}}$ the noise matrix. Let $\mathbf{r}(\boldsymbol{\theta})$ be the vectorised residual
\[
    \mathbf{r}(\boldsymbol{\theta}) = \mathrm{vec}\left(\mathbf{X} - \widehat{\mathcal{I}}_{\mathbf{J}}(\boldsymbol{\theta}) \right) \in \mathbb{R}^{M(n+1)}.
\]
Note that 
\[
    \mathbf{r}(\boldsymbol{\theta}^{*}) = \mathrm{vec} \left( \mathcal{I}(\boldsymbol{\theta}^{*}) + \mathbf{E}_{\text{noise}} - \widehat{\mathcal{I}}_{\mathbf{J}}(\boldsymbol{\theta}^{*}) \right) = \mathrm{vec} \left( \mathbf{E}_{\text{noise}} - \mathbf{E}_{\text{coll}}(\boldsymbol{\theta}^{*}) \right).
\]

Consider now the differentiable MSE loss function
\[
    \mathcal{L}(\boldsymbol{\theta}) = \frac{1}{M(n+1)} \left\| \mathbf{X} - \widehat{\mathcal{I}}_{\mathbf{J}}(\boldsymbol{\theta}) \right\|_F^2  = \frac{1}{M(n+1)} \|\mathbf{r}(\boldsymbol{\theta})\|_{2}^{2}  = \frac{1}{M(n+1)} \mathbf{r}(\boldsymbol{\theta})^{\mathrm{T}} \mathbf{r}(\boldsymbol{\theta})
\]
where $\cdot^{\mathrm{T}}$ denotes transpose. During training, we seek parameters $\hat{\boldsymbol{\theta}}$ such that
\[
    \hat{\boldsymbol{\theta}} = \arg\min_{\boldsymbol{\theta}} \mathcal{L}(\boldsymbol{\theta}).
\]
At such a local minimiser $\hat{\boldsymbol{\theta}}$, the gradient of the loss vanishes:
\[
    \nabla \mathcal{L}(\hat{\boldsymbol{\theta}}) = 0.
\]
Assuming $\hat{\boldsymbol{\theta}}$ is close to $\boldsymbol{\theta}^{*}$, we write the first-order Taylor expansion of $\nabla \mathcal{L}$ around $\boldsymbol{\theta}^{*}$:
\[
    \nabla \mathcal{L}(\hat{\boldsymbol{\theta}}) = \nabla \mathcal{L}(\boldsymbol{\theta}^{*}) + \nabla^2 \mathcal{L}(\boldsymbol{\theta}^{*})(\hat{\boldsymbol{\theta}} - \boldsymbol{\theta}^{*}) + \mathcal{O}(\|\hat{\boldsymbol{\theta}} - \boldsymbol{\theta}^{*}\|^2) = 0,
\]
and we have, to first order, the approximation
\[
    \mathbf{H}(\boldsymbol{\theta}^{*})(\hat{\boldsymbol{\theta}} - \boldsymbol{\theta}^{*}) \approx -\nabla \mathcal{L}(\boldsymbol{\theta}^{*}),
\]
where $\mathbf{H}(\boldsymbol{\theta}) = \nabla^2 \mathcal{L}(\boldsymbol{\theta})  \in \mathbb{R}^{d \times d}$ is the Hessian matrix of $\mathcal{L}(\boldsymbol{\theta})$. From this, we can estimate the parameter error:
\[
    \hat{\boldsymbol{\theta}} - \boldsymbol{\theta}^{*} \approx -\mathbf{H}^{\dagger}(\boldsymbol{\theta}^{*}) \  \nabla \mathcal{L}(\boldsymbol{\theta}^{*}),
\]
where $\cdot^{\dagger}$ denotes the Moore--Penrose pseudoinverse. 

Taking the 2-norm and applying the submultiplicative property yields the approximate bound
\[
    \| \hat{\boldsymbol{\theta}} - \boldsymbol{\theta}^{*}\|_{2} \lessapprox \| 
    \mathbf{\mathbf{H}}^{\dagger}(\boldsymbol{\theta}^{*}) \|_{2} \| \nabla \mathcal{L}(\boldsymbol{\theta}^{*}) \|_{2} = \frac{\| \nabla \mathcal{L}(\boldsymbol{\theta}^{*})
    \|_{2}}{\sigma_{\min}(\mathbf{H}(\boldsymbol{\theta}^{*}))},
\]
where $\sigma_{\mathrm{min}}(\cdot)$ denotes smallest singular value. This approximate bound is locally valid under the assumption that $\hat{\boldsymbol{\theta}}$ lies in a neighbourhood of $\boldsymbol{\theta}^{*}$ where the quadratic approximation of $\mathcal{L}$ is accurate. 

We now write $\nabla \mathcal{L}(\boldsymbol{\theta}^{*})$ in terms of the collocation error $\mathbf{E}_{\text{coll}}(\boldsymbol{\theta})$. The gradient of the loss is
\[
    \nabla \mathcal{L}(\boldsymbol{\theta}) = \frac{\partial}{\partial \boldsymbol{\theta}} \left[ \frac{1}{M(n+1)} \mathbf{r}(\boldsymbol{\theta})^{\mathrm{T}} \mathbf{r}(\boldsymbol{\theta}) \right] = A(M,n) \ \mathbf{F}(\boldsymbol{\theta})^{\mathrm{T}} \mathbf{r}(\boldsymbol{\theta}) \in \mathbb{R}^{d}
\]
where
\[
    \mathbf{F}(\boldsymbol{\theta}) = \frac{\partial \mathbf{r}(\boldsymbol{\theta})}{\partial \boldsymbol{\theta}} \in \mathbb{R}^{M(n+1) \times d},
\]
is the Jacobian of the forward map, measuring the sensitivity of observations to parameter changes, and with $A(M,n) = \frac{2}{M(n+1)}$ a constant depending on the size of the observations matrix. 

 Evaluating the gradient of the loss function at the ground-truth $\boldsymbol{\theta}^{*}$,
\[
    \nabla \mathcal{L}(\boldsymbol{\theta}^{*}) =  A(M,n) \ \mathbf{F}(\boldsymbol{\theta}^{*})^{\mathrm{T}} \mathrm{vec} \left( \mathbf{E}_{\text{noise}} - \mathbf{E}_{\text{coll}}(\boldsymbol{\theta}^{*}) \right),
\]
taking the 2-norm of this,
\begin{align*}
    \| \nabla \mathcal{L}(\boldsymbol{\theta}^{*}) \|_{2}
    &\le A(M,n) \ \|\mathbf{F}(\boldsymbol{\theta}^{*})^{\mathrm{T}}\|_{2} \, \| \operatorname{vec} \left( \mathbf{E}_{\text{noise}} - \mathbf{E}_{\text{coll}}(\boldsymbol{\theta}^{*}) \right) \|_{2} \\
    &\le A(M,n) \ \|\mathbf{F}(\boldsymbol{\theta}^{*})\|_{2} \, \left( \| \mathbf{E}_{\text{noise}} \|_{F} + \| \mathbf{E}_{\text{coll}}(\boldsymbol{\theta}^{*})\|_{F} \right),
\end{align*}
and using the previously derived bound~\eqref{eq:coll-bound} on $ \| \mathbf{E}_{\text{coll}}(\boldsymbol{\theta}^{*})\|_{F}$, together with the noise matrix bound
\[
    \|\mathbf{E}_{\text{noise}}\|_F \le \varepsilon \, \sqrt{M (n+1)},
\]
produces
\begin{align*}
    \| \nabla \mathcal{L}(\boldsymbol{\theta}^{*}) \|_{2}
    &\leq A(M,n) \ \sigma_{\max}(\mathbf{F}(\boldsymbol{\theta}^{*}))
    \left[ \varepsilon\sqrt{M (n+1)} + \sqrt{M(n+1)} \left( \frac{\kappa_{\text{int}}}{n^4} + \varepsilon n^{1/2} \right) \right] \nonumber\\[4pt]
    &= \sigma_{\max}(\mathbf{F}(\boldsymbol{\theta}^{*}))
    \left[ \frac{2}{\sqrt{M(n+1)}} \left( \frac{\kappa_{\text{int}}}{n^4} + \varepsilon ( 1+n^{1/2}) \right) \right]. \label{eq:grad-bound}
\end{align*}
And we have that
\begin{equation} \label{eq:param-bound}
    \| \hat{\boldsymbol{\theta}} - \boldsymbol{\theta}^{*}\|_{2} \lessapprox \frac{\sigma_{\max}(\mathbf{F}(\boldsymbol{\theta}^{*}))}{\sigma_{\min}(\mathbf{H}(\boldsymbol{\theta}^{*}))}
    \left[ \frac{2}{\sqrt{M(n+1)}} \left( \frac{\kappa_{\text{int}}}{n^4} + \varepsilon ( 1+n^{1/2}) \right) \right].
\end{equation}
This approximate bound separates the recovery error into two factors: a data-error factor
\[
    \left[ \frac{2}{\sqrt{M(n+1)}} \left( \frac{\kappa_{\text{int}}}{n^4} + \varepsilon ( 1+n^{1/2}) \right) \right]
\]
determined by the noise level and the accuracy of the collocation approximation; and a term
\[
    \frac{\sigma_{\max}(\mathbf{F}(\boldsymbol{\theta}^{*}))}{\sigma_{\min}(\mathbf{H}(\boldsymbol{\theta}^{*}))}
\]
depending on the conditioning of the Hessian matrix.

In the following subsection, we analyse how this conditioning is sensitive to particular structural properties of some CRN mechanisms. 

\subsection{Structural sources of ill-conditioning}
\label{sec:ill-conditioning}
The Hessian matrix
\[
    \mathbf{H}(\boldsymbol{\theta}) = \frac{\partial}{\partial \boldsymbol{\theta}} \left[ A(M,n) \ \mathbf{F}(\boldsymbol{\theta})^{\mathrm{T}} \mathbf{r}(\boldsymbol{\theta}) \right] = A(M,n) \left[ \frac{\partial \mathbf{F}(\boldsymbol{\theta})^{\mathrm{T}}}{\partial \boldsymbol{\theta}} \mathbf{r}(\boldsymbol{\theta}) + \mathbf{F}(\boldsymbol{\theta})^{\mathrm{T}} \mathbf{F}(\boldsymbol{\theta}) \right],
\]
inherits its rank and conditioning from the CRN architecture through the Jacobian matrix \(\mathbf{F}(\boldsymbol{\theta})\). In particular, we analyse three kinds of CRN structural deficiencies that impact the rank of \(\mathbf{H}(\boldsymbol{\theta})\) in different ways:
\begin{itemize}
    \item \textbf{Conservation laws between species:}
    If the CRN to be discovered follows a conservation law
    \[
        \sum_{S_\beta \in \mathcal{S}_{\mathrm{moiety}}} \dot{x}_{\beta}(t) =
        0 \Rightarrow \sum_{S_\beta \in \mathcal{S}_{\mathrm{moiety}}} x_{\beta}(t) =
        \mathrm{const}, \quad t \geq 0, \quad \mathcal{S}_{\mathrm{moiety}} \subseteq \mathcal{S},
    \]
    then we can find a nonzero vector $\mathbf{c}$ such that
    \[
        \mathbf c^\mathrm{T} \dot{\mathbf x}(t)=0
    \]
    for all $t$. Using the CRNN dynamical equation
    \(
        \dot{\mathbf x}(t)= \mathbf{Q}^{(\Delta)}\,\mathbf{R}(\mathbf{X}),
    \)
    we can see that this implies
    \[
        \mathbf c^\mathrm{T}\mathbf Q^{(\Delta)}=\mathbf 0.
    \]
    Thus, if our CRN has conservation laws between species, we will have both linear dependences among the features and among the parameters.

    There exists a constant left-null vector of \(\mathbf{F}(\boldsymbol{\theta})\),
    \[
        \mathbf{c}_l^{\mathrm{T}} \mathbf{F}(\boldsymbol{\theta}) = \mathbf{0}, \qquad \mathbf{c}_l \neq \mathbf{0}, \qquad \mathbf{c}_l \in \mathbb{R}^{M(n+1)},
    \]
    that encodes the linear dependence over observations at the vectorised level (rows of \(\mathbf{F}\)). This vector satisfies
    \[
        \mathbf{F}(\boldsymbol{\theta}) \mathbf{F}(\boldsymbol{\theta})^{\mathrm{T}} \mathbf{c}_l = \mathbf{F}(\boldsymbol{\theta}) \bigl(\mathbf{c}_l^{\mathrm{T}} \mathbf{F}(\boldsymbol{\theta})\bigr)^{\mathrm{T}} = \mathbf{0},
    \]
    so
    \[
        \mathrm{rank}(\mathbf{F}(\boldsymbol{\theta})\mathbf{F}(\boldsymbol{\theta})^{\mathrm{T}}) = \mathrm{rank}(\mathbf{F}(\boldsymbol{\theta})) < M(n+1).
    \]
    The Hessian matrix $\mathbf{H}(\boldsymbol{\theta})$ inherits this rank deficiency via
    \[
        \mathrm{rank}(\mathbf{H}(\boldsymbol{\theta})) \leq \mathrm{rank}(\mathbf{F}(\boldsymbol{\theta})) < M(n+1).
    \]

    There also exists a constant right-null vector of \(\mathbf{F}(\boldsymbol{\theta})\),
    \[
        \mathbf{F}(\boldsymbol{\theta}) \mathbf{c}_r = \mathbf{0}, \qquad \mathbf{c}_r \neq \mathbf{0},\qquad \mathbf{c}_r \in \mathbb{R}^{d},
    \]
    that encodes the linear dependence over parameters at the vectorised level (columns of \(\mathbf{F}\)).
    
    Because this null vector is constant, we have that
    \[
        \frac{\partial \mathbf{F}(\boldsymbol{\theta})}{\partial \boldsymbol{\theta}} \mathbf{c}_r = \frac{\partial}{\partial \boldsymbol{\theta}} \left(\mathbf{F} (\boldsymbol{\theta})\mathbf{c}_r\right) = \mathbf{0}.
    \]

    Furthermore, \(\mathbf{c}_r\) is a left-null vector of \(\mathbf{H}(\boldsymbol{\theta})\),
    \[
        \mathbf{c}_r^{\mathrm{T}} \mathbf{H}(\boldsymbol{\theta}) = \mathbf{0},
    \]
    because
    \[
        \mathbf{c}_r^{\mathrm{T}} \mathbf{F}(\boldsymbol{\theta})^{\mathrm{T}} \mathbf{F}(\boldsymbol{\theta}) = \bigl(\mathbf{F}(\boldsymbol{\theta})\,\mathbf{c}_r\bigr)^{\mathrm{T}} \mathbf{F}(\boldsymbol{\theta}) = \mathbf{0},
    \]
    and
    \[
        \mathbf{c}_r^{\mathrm{T}} \frac{\partial \mathbf{F}(\boldsymbol{\theta})^{\mathrm{T}}}{\partial \boldsymbol{\theta}} \mathbf{r}(\boldsymbol{\theta}) = \mathbf{r}(\boldsymbol{\theta})^{\mathrm{T}} \frac{\partial \mathbf{F}(\boldsymbol{\theta})}{\partial \boldsymbol{\theta}} \mathbf{c}_r = \mathbf{r}(\boldsymbol{\theta})^{\mathrm{T}} \frac{\partial}{\partial \boldsymbol{\theta}} \bigl(\mathbf{F}(\boldsymbol{\theta})\,\mathbf{c}_r\bigr) = 0.
    \]
    Since \(\mathbf{c}_r \neq \mathbf{0}\), this is a nontrivial null vector of \(\mathbf{H}(\boldsymbol{\theta})\), so when our CRN has linear dependences between parameters, the Hessian matrix inherits this as
    \[
        \mathrm{rank}(\mathbf{H}(\boldsymbol{\theta}))< d.
    \]

    When the CRN has linear dependences both between the parameters and between the concentrations, we can thus conclude that
    \[
        \mathrm{rank}(\mathbf{H}(\boldsymbol{\theta})) < \min(d, M(n+1)),
    \]

    \item \textbf{Presence of reversible reactions.}
    Consider a pair of reversible reactions
    \[
        (q_i,q_j), \ (q_j,q_i) \in \mathcal{R}
    \]
    with forward and backward rate constants \(k_{i,j}\) and \(k_{j,i}\), respectively, where \((q_i,q_j)\) denotes the reaction converting complex \(q_i\) into complex \(q_j\). Their contribution to the dynamics is
    \[
        \mathbf Q^{(\Delta)}_{:,(i,j)}
        k_{i,j}d_i(\mathbf x)
        +
        \mathbf Q^{(\Delta)}_{:,(j,i)}
        k_{j,i}d_j(\mathbf x).
    \]
    Since
    \[
        \mathbf Q^{(\Delta)}_{:,(j,i)} = -\mathbf Q^{(\Delta)}_{:,(i,j)},
    \]
    the two reactions contribute along the same stoichiometric direction with opposite signs, introducing a linear dependence between the columns of $\mathbf Q^{(\Delta)}$.

    There exists a right-null vector of $\mathbf Q^{(\Delta)}$,
    \[
        \mathbf Q^{(\Delta)} \mathbf{c}' = 0, \qquad \mathbf{c}' \neq 0, \qquad \mathbf{c}' \in \mathbb{R}^{R}.
    \]
    Which implies that there also exists a constant right-null vector of \(\mathbf{F}(\boldsymbol{\theta})\),
    \[
        \mathbf{F}(\boldsymbol{\theta}) \mathbf{c}'_r = \mathbf{0}, \qquad \mathbf{c}'_r \neq \mathbf{0},\qquad \mathbf{c}'_r \in \mathbb{R}^{d}.
    \]
    As in the previous case, the impact of linear dependences on the parameters over the Hessian rank is 
    \[
        \mathrm{rank}(\mathbf{H}(\boldsymbol{\theta}))< d.
    \]

    \item \textbf{Shared reactants.} 
    When multiple reactions share the same reactant complex, e.g.,
    \[
        (q_i, q_{j_1}), \ (q_i, q_{j_2}) \in \mathcal{R},
    \]
    with rate constants \(k_{i,j_1}\) and \(k_{i,j_2}\) respectively, these reactions depend on the same mass-action monomial
    \[
        d_i(\mathbf x)
        =
        \prod_{\alpha=1}^{M}
        x_\alpha^{Q_{\alpha,i}}.
    \]
    Their combined contribution is
    \[
        \left(
            \mathbf Q^{(\Delta)}_{:,(i,j_1)}k_{i,j_1}
            +
            \mathbf Q^{(\Delta)}_{:,(i,j_2)}k_{i,j_2}
        \right)
        d_i(\mathbf x).
    \]
    Hence we may detect only the resulting linear combination of the two reaction contributions, rather than each reaction individually. Unlike the two previous cases, this does not produce an exact null vector of $\mathbf F$; instead it produces an alternative representation of the dynamics in which the number of reactions is reduced and the corresponding stoichiometric columns are merged.
\end{itemize}

\subsection{An illustrative example}

\begin{figure}[H]
    \centering
    \begin{tikzpicture}
        \node (A) at (0,0) {$\sA+\scat$};
        \node (B) at (3,0) {$\scatA$};
        \node (C) at (6,0) {$\sP+\scat$};
        \draw[->] (A) to[out=45,in=135]  node[above] {$k_{1}$} (B);
        \draw[<-] (A) to[out=-45,in=-135]  node[below] {$k_{-1}$} (B);
        \draw[->] (B) to[out=45,in=135]   node[above] {$k_{2}$} (C);
        \draw[<-] (B) to[out=-45,in=-135]   node[below] {$k_{-2}$} (C);
    \end{tikzpicture}
    \caption{Graph of the fully-reversible Michaelis-Menten mechanism}
    \label{fig:M1-Graph}
\end{figure}
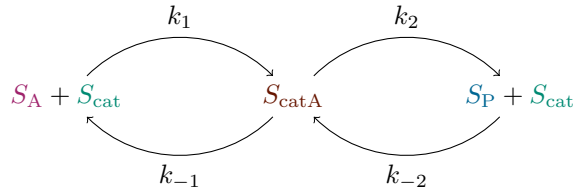

We now exemplify the aforementioned potential sources of Hessian ill-conditioning induced by CRN architecture by considering the fully-reversible Michaelis-Menten reaction network, or M1 model \cite{bures2023organic}, shown in Figure~\ref{fig:M1-Graph}.

\noindent The log-exp reparametrised reaction-explicit ODE system~\eqref{eq:NODE} for the M1 model is
\begin{equation} \label{eq:M1-react-ODE}
    \overbrace{\begin{bmatrix}
        \dxA     \\
        \dxP     \\
        \dxcat   \\
        \dxcatA
    \end{bmatrix}}^{\mathbf{\dot{x}}(t)} =
    \overbrace{\begin{bmatrix}
        -1 & +1 &  0 &  0 \\
         0 &  0 & +1 & -1 \\
        -1 & +1 & +1 & -1 \\
        +1 & -1 & -1 & +1
    \end{bmatrix}}^{\mathbf{Q}^{(\Delta)}}
    \exp \left(
    \overbrace{\begin{bmatrix}
        1 & 0 & 1 & 0 \\
        0 & 0 & 0 & 1 \\
        0 & 0 & 0 & 1 \\
        0 & 1 & 1 & 0
    \end{bmatrix}}^{(\mathbf{Q}^{(\mathrm{Re})})^{T}} \overbrace{\begin{bmatrix}
        \log(\xA)     \\
        \log(\xP)      \\
        \log(\xcat)    \\
        \log(\xcatA)
    \end{bmatrix}}^{\mathbf{{x}}(t)}   +
    \overbrace{\begin{bmatrix}
        \log(k_{1}) \\ \log(k_{-1}) \\ \log(k_{2}) \\ \log(k_{-2})
    \end{bmatrix}}^{\mathbf{b}}
    \right).
\end{equation}

Some examples of structural ill-conditioning of~\eqref{eq:M1-react-ODE} that are inherited by the Hessian matrix are the following.

\paragraph{Conservation laws between species.} 
In the case of~\eqref{eq:M1-react-ODE}, we have two moieties:
\begin{align*}
    \mathcal{S}_{\mathrm{1}} = \{\scat,\scatA\} \ \ &\mathrm{s.t.} \ \  \dxcat+\dxcatA = 0;  \\
    \mathcal{S}_{\mathrm{2}} = \{\sA, \sP, \scat\} \ \ &\mathrm{s.t.} \ \  \dxA+\dxP+\dxcat = 0,
\end{align*}
which in turn introduce the linear dependencies
\[
    \mathbf{Q}^{(\Delta)}_{3,:} = - \mathbf{Q}^{(\Delta)}_{4,:},
    \qquad
    \mathbf{Q}^{(\Delta)}_{1,:} + \mathbf{Q}^{(\Delta)}_{2,:} = -\mathbf{Q}^{(\Delta)}_{3,:}
\]
between the rows of $\mathbf{Q}^{(\Delta)}$, and
\[
    \mathbf{Q}^{(\mathrm{Re})}_{1,:} + \mathbf{Q}^{(\mathrm{Re})}_{2,:} = \mathbf{Q}^{(\mathrm{Re})}_{3,:}
\]
between the rows of $\mathbf{Q}^{(\mathrm{Re})}$.

\paragraph{Reversible reactions} 
In the case of~\eqref{eq:M1-react-ODE}, we have two pairs of reversible reactions,
\[
    \bigl((\sA+\scat,\scatA), (\scatA,\sA+\scat)\bigr),
    \qquad
    \bigl((\scatA,\sP+\scat), (\sP+\scat,\scatA)\bigr),
\]
which in turn imply the equalities
\[
    \mathbf{Q}^{(\Delta)}_{:,1} = - \mathbf{Q}^{(\Delta)}_{:,2},
    \qquad
    \mathbf{Q}^{(\Delta)}_{:,3} = - \mathbf{Q}^{(\Delta)}_{:,4}
\]
between the columns of $\mathbf{Q}^{(\Delta)}$.

\paragraph{Shared reactants.}
In the case of~\eqref{eq:M1-react-ODE}, the two reactions $(\scatA,\sA+\scat)$ and $(\scatA,\sP+\scat)$ share the reactant \(\scatA\). This allows for alternative representations of the ODE system~\eqref{eq:M1-react-ODE}. For instance, merging the two reactions out of $\scatA$ into a single effective reaction gives the three-reaction system
\begin{equation*}
    \begin{bmatrix}
        \dxA     \\
        \dxP     \\
        \dxcat   \\
        \dxcatA
    \end{bmatrix} =
    \begin{bmatrix}
        -1 & +c_{1} &  0 \\
         0 & +c_{2} & -1 \\
        -1 & +1 & -1 \\
        +1 & -1 & +1
    \end{bmatrix}
    \exp \left(
    \begin{bmatrix}
        1 & 0 & 1 & 0 \\
        0 & 0 & 0 & 1 \\
        0 & 1 & 1 & 0
    \end{bmatrix}
    \begin{bmatrix}
        \log(\xA)     \\
        \log(\xP)      \\
        \log(\xcat)    \\
        \log(\xcatA)
    \end{bmatrix}  +
    \begin{bmatrix}
        \log(k_{1}) \\ \log(k_{-1}+k_{2}) \\ \log(k_{-2})
    \end{bmatrix}
    \right),
\end{equation*}
with
\[
    c_{1} = \frac{k_{-1}}{k_{-1}+k_{2}}, \quad c_{2} = \frac{k_{2}}{k_{-1}+k_{2}},
\]
which reproduces the dynamics~\eqref{eq:M1-react-ODE} exactly but has a different stoichiometric matrix: the third column is now a linear combination of the two original columns, and $c_1, c_2$ are non-integer in general. Thus, it can happen that we learn an ODE that confidently fits the training data and produces reliable predictions of the system's evolution, while the true network topology remains obscured: because $c_{1}, c_{2}$ are non-integer, the resulting stoichiometric matrix has no valid chemical interpretation, and the original four-reaction network is not identified. This is a qualitatively different problem than the two preceding cases.



\subsection{Reduction of parameter space}
CRN structure can also aid in simplifying the inverse problem of our interest. For instance, if we assume that the CRN to be discovered does not include autocatalytic reactions, i.e., reactions $(q_i,q_j)$ such that $Q_{\alpha,i} \neq 0$ and $Q_{\alpha,j} \neq 0$ for some $S_\alpha \in \mathcal{S}$, e.g.
\[
    S_{1}+S_{2} \rightleftarrows 2S_{2} \quad \mathrm{is \ an \ autocatalytic \ reaction}.
\]
Under this assumption, the net stoichiometry $\mathbf{Q}^{(\Delta)}$ matrix uniquely identifies which species are consumed (negative entries) and which are produced (positive entries) and we can therefore directly relate $\mathbf{Q}^{(\mathrm{Re})}$ and $\mathbf{Q}^{(\Delta)}$ through the identity
\[
    \mathbf{Q}^{(\mathrm{Re})} = \max \left(-\mathbf{Q}^{(\Delta)},\, 0\right),
\]
which extracts the consumed species from the net stoichiometry and sets them as reactants. For these systems, the learnable parameter space reduces to
\[
    \boldsymbol{\theta} \rightarrow \boldsymbol{\theta}' = 
    (\mathbf{Q}^{(\Delta)},\, \mathbf{b}) \in \mathbb{R}^{R(M+1)}.
\]

\section{Numerical experiments}\label{sec:experiments}

We now validate our iCRNN approach through a series of numerical experiments. We compare the performance of training a CRNN using an iterative ODE solver against using the integration matrix~$\mathbf{J}$, and assess recovery accuracy, training cost, and robustness to measurement noise on two benchmark mechanisms.

Our numerical experiments were performed in Julia. We use the original CRNN approach of Ji \& Deng as a benchmark for comparison, with our implementation closely following that of the author's available at \url{https://github.com/DENG-MIT/CRNN}. 
Our Julia repository is available at 
\begin{center}
    \url{https://github.com/nla-group/iCRNN}
\end{center}
The  runtimes correspond to running the code on a MacBook Air M2 laptop.

\subsection{M1 model}

We generate $100$ independent batches of synthetic data for this system, each with a distinct set of kinetic constants. Each batch consists of $32$ experiments, divided into $24$ training experiments and $8$ validation experiments. Initial concentrations are assigned in two groups: 
\begin{itemize}
    \item for three quarters of the experiments we set $\xA(0) = 1.0$, $\xP(0) = \xcatA(0) = 0.0$, and draw $\xcat(0) \sim \mathcal{U}(0.1, 0.25)$;

    \item for the remaining quarter we draw $\xA(0) \sim \mathcal{U}(0.4, 0.8)$ and $\xcat(0) \sim \mathcal{U}(0.1, 0.25)$, set $\xP(0) = 1 - \xA(0)$, and set $\xcatA(0) = 0.0$.
\end{itemize}  
The kinetic constants are sampled independently from log-uniform distributions: 
\[
    k_1 \sim \log \mathcal{U}(10^{-1}, 10^{2}), \quad k_{-1} \sim \log \mathcal{U}(10^{-1}, 10^{0}), \quad k_2 \sim \log \mathcal{U}(10^{-1}, 10^{2}), \quad k_{-2} \sim \log \mathcal{U}(10^{-1}, 10^{0}).
\]
The training data is generated through numerical integration with LSODA \cite{petzold1983automatic} over a per-batch time interval $[0, t_{\max}]$ using $100$ uniformly spaced time points, where $t_{\max}$ is chosen individually for each batch to capture the full transient behaviour of the trajectories. Three noise levels are generated by multiplicative Gaussian perturbations of $0\%$, $2.5\%$, and $5\%$, yielding noiseless, low-noise, and high-noise datasets, respectively.

\begin{figure}[H]
    \centering
    \includegraphics[width=\linewidth]{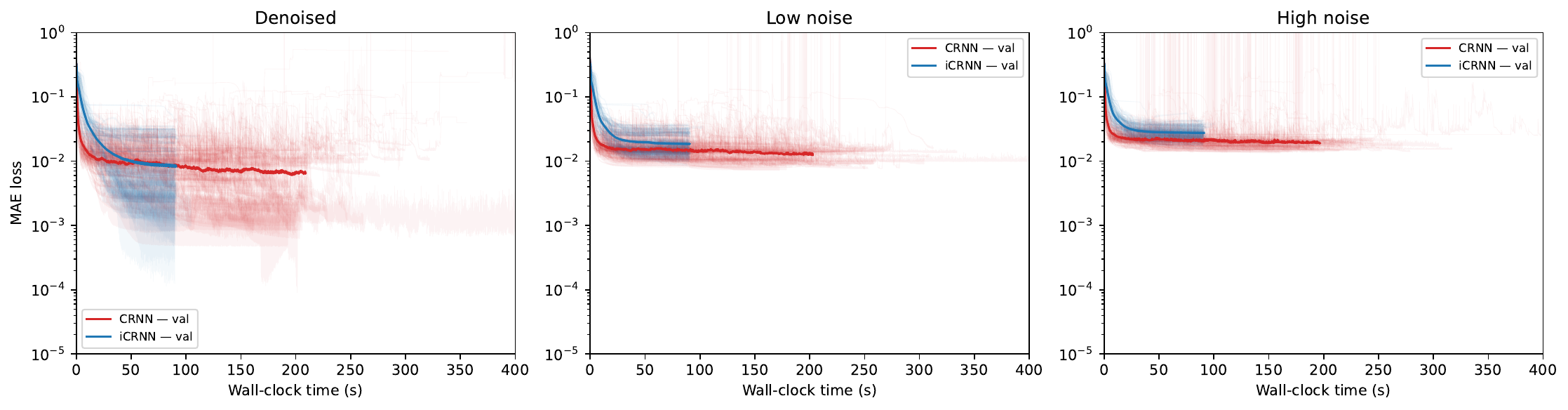}
    \caption{Validation loss versus wall-clock time for the M1 mechanism under MAE training, shown separately for the denoised, low-noise, and high-noise datasets. Solid lines are the geometric mean over $100$ kinetic parameter sets; faint lines are individual runs. iCRNN reaches its final loss in roughly half the wall-clock time of CRNN at every noise level, and the training loss history is significantly smoother.}
    \label{fig:M1-MAE-vs-time}
\end{figure}

\begin{figure}[H]
    \centering
    \includegraphics[width=\linewidth]{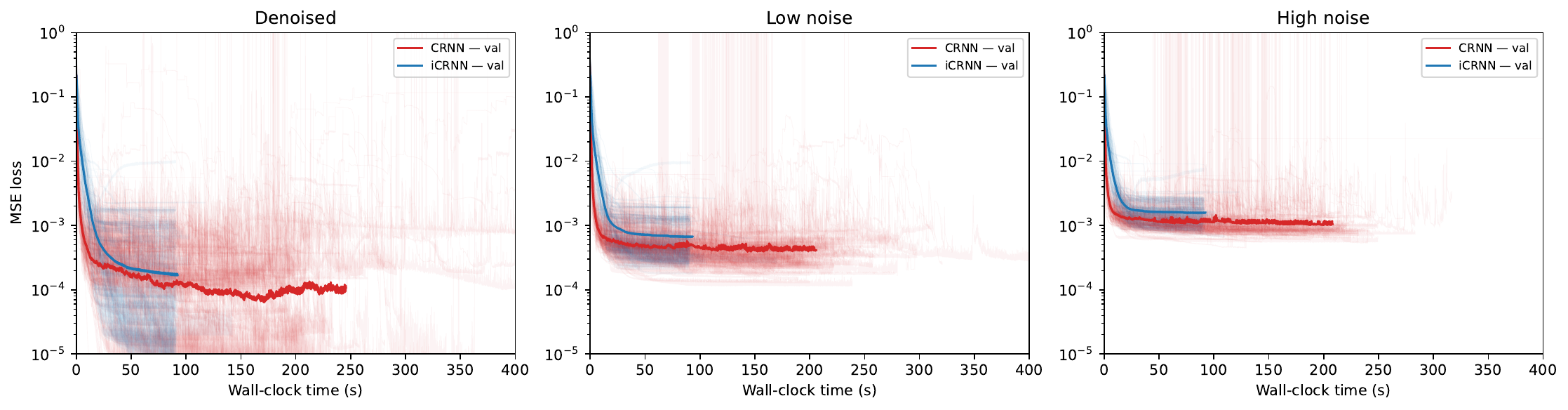}
    \caption{As Figure~\ref{fig:M1-MAE-vs-time}, but using MSE as the loss when training. There is very little difference between MAE and MSE in the qualitative training behaviour.}
    \label{fig:M1-MSE-vs-time}
\end{figure}

All training runs of CRNN and iCRNN are performed over 20,000 epochs, in all cases performed in Julia with the same loss functions (MAE or MSE) and the same optimizer (AdamW). The results are shown in Figures~\ref{fig:M1-MAE-vs-time} and~\ref{fig:M1-MSE-vs-time}. 
There is a significant difference in the qualitative training behaviour of the two methods. The iCRNN loss decreases smoothly throughout training, and most individual runs remain tightly clustered around the mean. The CRNN, by contrast, exhibits pronounced run-to-run variability: individual trajectories fluctuate, sometimes by more than an order of magnitude between consecutive checkpoints. 

The observed qualitative difference reflects the different forward pass approach. iCRNN evaluates the reaction rates directly at the observed time points and integrates them through a fixed linear operator. Hence, the loss landscape is a smooth algebraic function of the parameters, and gradient updates produce correspondingly smooth changes in the loss. The ODE-solver approach instead differentiates through an adaptive integrator (Tsit5~\cite{tsitouras2011runge}) whose step sizes, tolerances, and error-control decisions depend on the current parameter values, introducing discontinuities and non-smoothness into the effective loss surface. Some runs of CRNN may also fail to complete, typically because the forward pass blows up. In summary, a single training run of iCRNN is usually more predictable than that of CRNN.

The smoothness and predictability of the iCRNN loss comes at a cost. Since the collocation operator is fixed and independent of the parameters, the forward pass is less flexible than an adaptive solver, and the CRNN approach consequently attains a slightly lower final loss at every noise level (Table~\ref{tab:M1-summary}). The two methods are therefore not in direct competition on loss alone: CRNN fits the training data more tightly, at roughly twice the computational cost, while iCRNN trades a small amount of fit for a smoother optimisation landscape, improved runtime, and complete convergence across all runs.

\begin{table}[H]
    \centering
    \arrayrulecolor{gray!60}

    \begin{tabular}{|c|c|c|c|c|c|c|}
        \hline
        \rowcolor{gray!15}
        \textbf{Method} &
        \textbf{Loss} &
        \textbf{Noise} &
        \textbf{Train loss} &
        \textbf{Val loss} &
        \textbf{Runtime (s)} &
        \textbf{Completed runs} \\
        \hline

        \multirow{6}{*}{CRNN}
        & \multirow{3}{*}{MAE}
        & Denoised   & $6.3\times 10^{-3}$  & $6.5\times 10^{-3}$  & 213.2 & 98/100 \\ \cline{3-7}
        & & Low noise  & $1.2\times 10^{-2}$  & $1.3\times 10^{-2}$  & 202.6 & 100/100 \\ \cline{3-7}
        & & High noise & $1.8\times 10^{-2}$  & $1.9\times 10^{-2}$  & 196.5 & 100/100 \\ \cline{2-7}

        & \multirow{3}{*}{MSE}
        & Denoised   & $8.7\times 10^{-5}$  & $1.0\times 10^{-4}$  & 241.6 & 96/100 \\ \cline{3-7}
        & & Low noise  & $3.6\times 10^{-4}$  & $4.1\times 10^{-4}$  & 205.7 & 99/100 \\ \cline{3-7}
        & & High noise & $9.4\times 10^{-4}$  & $1.1\times 10^{-3}$  & 208.1 & 100/100 \\
        \hline

        \multirow{6}{*}{iCRNN}
        & \multirow{3}{*}{MAE}
        & Denoised   & $7.7\times 10^{-3}$  & $8.3\times 10^{-3}$  & 90.9  & 100/100 \\ \cline{3-7}
        & & Low noise  & $1.7\times 10^{-2}$  & $1.9\times 10^{-2}$  & 90.5  & 100/100 \\ \cline{3-7}
        & & High noise & $2.6\times 10^{-2}$  & $2.8\times 10^{-2}$  & 90.8  & 100/100 \\ \cline{2-7}

        & \multirow{3}{*}{MSE}
        & Denoised   & $1.5\times 10^{-4}$  & $1.7\times 10^{-4}$  & 92.2  & 100/100 \\ \cline{3-7}
        & & Low noise  & $5.8\times 10^{-4}$  & $6.7\times 10^{-4}$  & 93.2  & 100/100 \\ \cline{3-7}
        & & High noise & $1.4\times 10^{-3}$  & $1.6\times 10^{-3}$  & 92.2  & 100/100 \\
        \hline
    \end{tabular}
    \caption{Final training and validation losses (geometric mean over $100$ kinetic parameter sets) and mean wall-clock runtime for the M1 mechanism, across both training losses and all three noise levels. Runtime is reported as an arithmetic mean. ``Completed runs'' gives the number of runs that finished and produced a finite final loss, out of $100$.}
    \label{tab:M1-summary}
\end{table}





\subsection{Example model (EM)}
The second benchmark is the five-species example model (EM) of Ji \& Deng~\cite{ji2021autonomous}, comprising the reactions
\begin{align*}
    2A &\xrightarrow{k_1} B ; \\
    A &\xrightarrow{k_2} C ; \\
    C &\xrightarrow{k_3} D ; \\
    B+D &\xrightarrow{k_4} E.
\end{align*}

We generate $50$ independent batches of synthetic data, each with a distinct set of kinetic constants and $30$ experiments per batch,
divided into $20$ training experiments and $10$ validation experiments. Initial concentrations are drawn independently for each experiment as $x_A(0) \sim \mathcal{U}(0.2, 1.2)$ and $x_B(0) \sim \mathcal{U}(0.2, 1.2)$, with $x_C(0) = x_D(0) = x_E(0) = 0.0$. 

The kinetic constants $[k_1, k_2, k_3, k_4]$ are sampled independently from log-uniform distributions over $[10^{-1}, 10^{1}]$. 

The training data is generated through numerical integration with LSODA \cite{petzold1983automatic} over a per-batch time interval $[0, t_{\max}]$ using $100$ uniformly spaced time points, with $t_{\max}$ chosen individually for each batch to capture the full transient behaviour of the trajectories. As for M1, three noise levels are generated by multiplicative Gaussian perturbations of $0\%$, $2.5\%$, and $5\%$.

\begin{figure}[H]
    \centering
    \includegraphics[width=\linewidth]{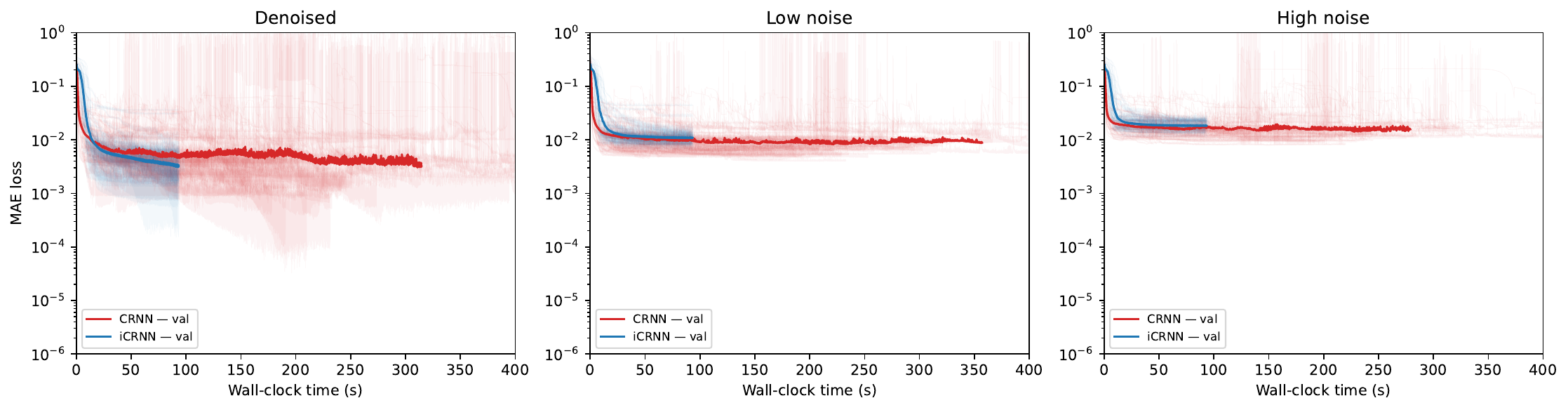}
    \caption{Validation loss versus wall-clock time for the EM mechanism under MAE training, shown separately for the denoised, low-noise, and high-noise datasets. Solid lines are the geometric mean over $50$ kinetic parameter sets; faint lines are individual runs. iCRNN reaches its final loss in roughly half the wall-clock time of CRNN at every noise level.}
    \label{fig:EM-MAE-vs-time}
\end{figure}

\begin{figure}[H]
    \centering
    \includegraphics[width=\linewidth]{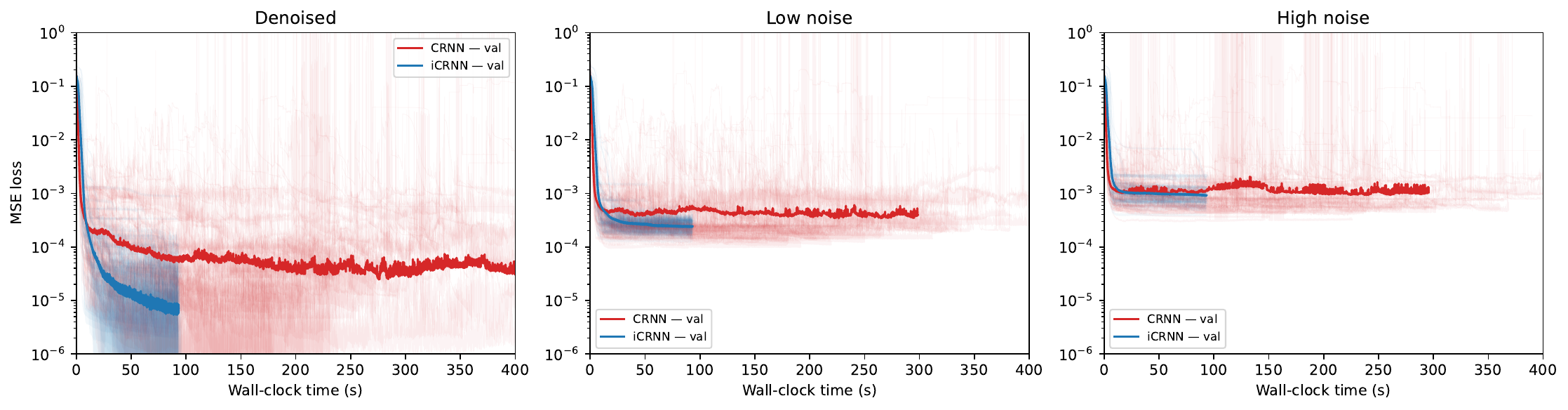}
    \caption{As Figure~\ref{fig:EM-MAE-vs-time}, but using MSE as the loss when training. There is very little difference between MAE and MSE in the qualitative training behaviour.}
    \label{fig:EM-MSE-vs-time}
\end{figure}

\begin{table}[H]
    \centering
    \arrayrulecolor{gray!60}

    \begin{tabular}{|c|c|c|c|c|c|c|}
        \hline
        \rowcolor{gray!15}
        \textbf{Method} &
        \textbf{Loss} &
        \textbf{Noise} &
        \textbf{Train loss} &
        \textbf{Val loss} &
        \textbf{Runtime (s)} &
        \textbf{Completed runs} \\
        \hline

        \multirow{6}{*}{CRNN}
        & \multirow{3}{*}{MAE}
        & Denoised   & $3.2\times 10^{-3}$  & $3.2\times 10^{-3}$  & 321.0 & 49/50 \\ \cline{3-7}
        & & Low noise  & $8.8\times 10^{-3}$  & $8.9\times 10^{-3}$  & 357.0 & 50/50 \\ \cline{3-7}
        & & High noise & $1.5\times 10^{-2}$  & $1.6\times 10^{-2}$  & 284.6 & 49/50 \\ \cline{2-7}

        & \multirow{3}{*}{MSE}
        & Denoised   & $3.5\times 10^{-5}$  & $2.9\times 10^{-5}$  & 409.2 & 48/50 \\ \cline{3-7}
        & & Low noise  & $4.7\times 10^{-4}$  & $3.9\times 10^{-4}$  & 298.5 & 49/50 \\ \cline{3-7}
        & & High noise & $8.5\times 10^{-4}$  & $9.8\times 10^{-4}$  & 284.7 & 49/50 \\
        \hline

        \multirow{6}{*}{iCRNN}
        & \multirow{3}{*}{MAE}
        & Denoised   & $3.1\times 10^{-3}$  & $3.3\times 10^{-3}$  & 93.2  & 50/50 \\ \cline{3-7}
        & & Low noise  & $1.0\times 10^{-2}$  & $1.1\times 10^{-2}$  & 93.2  & 50/50 \\ \cline{3-7}
        & & High noise & $1.6\times 10^{-2}$  & $1.8\times 10^{-2}$  & 93.2  & 50/50 \\ \cline{2-7}

        & \multirow{3}{*}{MSE}
        & Denoised   & $6.2\times 10^{-6}$  & $7.3\times 10^{-6}$  & 93.3  & 50/50 \\ \cline{3-7}
        & & Low noise  & $2.2\times 10^{-4}$  & $2.4\times 10^{-4}$  & 93.1  & 50/50 \\ \cline{3-7}
        & & High noise & $7.3\times 10^{-4}$  & $9.1\times 10^{-4}$  & 93.2  & 50/50 \\
        \hline
    \end{tabular}
    \caption{Final training and validation losses (geometric mean over $50$
    kinetic parameter sets) and mean wall-clock runtime for the EM
    mechanism, across both training losses and all three noise levels.
    Runtime is reported as an arithmetic mean. ``Completed runs'' gives the number of runs that finished
    and produced a finite final loss, out of $50$ per condition.}
    \label{tab:EM-summary}
\end{table}






Table~\ref{tab:EM-summary} reports the final losses and runtimes for the EM mechanism. For this experiment, the runtime advantage of iCRNN over CRNN is more evident: iCRNN completes every training in approximately $93$~s, whereas CRNN requires between $285$ and $409$~s, a speedup of $3.1$ to $4.4\times$. This gap might indicate that the advantage of the algebraic forward pass of iCRNN grows as we consider CRNs with greater number of species. 

On MSE, iCRNN also attains lower final loss at every noise level, with the margin largest on the denoised data ($5.6\times$) and shrinking as noise increases ($2.1\times$ at low noise, $1.2\times$ at high noise). This is consistent with the data-error factor in the bound~\eqref{eq:param-bound}. On MAE, by contrast, the two methods are comparable, with iCRNN marginally better on denoised data and CRNN marginally better at low and high noise; the differences are within a factor of $1.2$ and are unlikely to be meaningful.

The qualitative smoothness of the iCRNN loss curves observed for M1 is reproduced for EM (Figures~\ref{fig:EM-MAE-vs-time} and~\ref{fig:EM-MSE-vs-time}), as is the completion advantage: iCRNN finished all $50$ runs at every noise level, whereas CRNN failed to converge on $1$ to $2$ runs per condition.

\section{Conclusions}\label{sec:conclusions}
We have introduced iCRNN, an integral reformulation of the chemical reaction neural network in which the forward pass is entirely algebraic. By evaluating the learned reaction rates at the observed time points and integrating them through a fixed cubic spline collocation operator, iCRNN removes the repeated adaptive ODE solves that dominate the cost of standard CRNN training. The resulting loss landscape is a smooth algebraic function of the parameters, so gradients are well-defined everywhere and do not depend on solver tolerances or stiffness.

Across two benchmark mechanisms with comparable implementations in Julia, iCRNN trains in roughly half the wall-clock time of CRNN on the four-species M1 mechanism and between three and four times faster on the five-species EM mechanism, with speedup seeming to grow as the system size increases. iCRNN also completed every training run on both mechanisms, whereas CRNN failed to converge on up to four runs per condition on M1 and up to two on EM. On final loss, the two methods are comparable: CRNN attains slightly lower final loss on M1, while iCRNN attains lower loss on EM trained with MSE and matches CRNN on EM trained with MAE. The overall picture is therefore a trade-off rather than a dominance: iCRNN exchanges a small amount of flexibility in fitting the training trajectory for a smoother optimisation landscape, substantially reduced runtime, and more reliable behaviour.

The error analysis clarifies when this trade-off is favourable and when it is not. We derived a bound on the parameter recovery error and separated it into a conditioning factor, determined by the structure of the network, and a data-error factor, determined by the noise level and the accuracy of the collocation approximation. This decomposition makes precise the intuition that chemical model recovery can be poor even when the trajectory predictor residual is small, and it identifies three structural sources of ill-conditioning. Conservation laws and reversible reactions introduce linear dependencies that inflate parameter uncertainty but leave the topology recoverable in principle. Shared reactants are qualitatively different: they allow the same trajectory to be reproduced by a network with a different stoichiometric structure, so the model may fit the data confidently while the true graph remains unidentifiable. This distinction is important in practice because it means that low training and validation loss are not by themselves evidence that the recovered mechanism is correct.

Several limitations remain. The iCRNN collocation operator is fixed and assumes uniformly spaced observations; extending the formulation to non-uniform or adaptive time grids would broaden its applicability. The error analysis is first-order and local, valid in a neighbourhood of the ground-truth parameters where the quadratic approximation of the loss is accurate; a Bayesian formulation of the recovery problem would provide sharper guarantees and a principled way to quantify uncertainty in the recovered rates and topology. The structural obstructions identified in this work are, by their nature, not resolved by a better optimiser or a faster forward pass. Shared reactants, in particular, produce alternative parametrisations that are indistinguishable from the data, and recovering the correct topology requires additional information, such as prior constraints on the admissible reaction structures. Extending the approach to settings where observed trajectories
are sparse, meaning that we only measure some species evolution, or irregularly timed, also remains open.

We emphasize again that our collocation formulation of the integral forward pass is not specific to mass-action kinetics or to chemical reaction networks, and may extend to other neural ODEs with the same fast, smooth, solver-free training regime. We leave this for future work.

\bibliographystyle{plain}
\bibliography{references}

\end{document}